\documentclass[12pt]{article}
\newtheorem{lemma}{Lemma}[section]
\newtheorem{theorem}[lemma]{Theorem}
\newtheorem{proposition}[lemma]{Proposition}

\newtheorem{definition}[lemma]{Definition}
\newenvironment{proof}{{\bf Proof}}{{\hfill $ \Box $}\vskip 4mm}
\newenvironment{remark}{\addtocounter{lemma}{1}
{\bf Remark \thelemma}}{{\hfill}\vskip 4mm}
\newenvironment{remarks}{\addtocounter{lemma}{1}
{\bf Remarks \thelemma}}{{\hfill}\vskip 4mm}

\newenvironment{examples}{\addtocounter{lemma}{1}
{\bf Examples \thelemma}}{{\hfill}\vskip 4mm}
\newcommand{\nc}{\newcommand}
\nc{\rnc}{\renewcommand}
\nc{\nt}{\newtheorem}

        {~\hfill~\fbox{}\end{trivlist}}

\nc{\thlabel}[1]{\label{theo:#1}}
\nc{\thref}[1]{Theorem~\ref{theo:#1}}
\nc{\selabel}[1]{\label{sect:#1}}
\nc{\seref}[1]{Section~\ref{sect:#1}}
\nc{\lelabel}[1]{\label{lemm:#1}}
\nc{\leref}[1]{Lemma~\ref{lemm:#1}}
\nc{\prlabel}[1]{\label{prop:#1}}
\nc{\prref}[1]{Proposition~\ref{prop:#1}}
\nc{\colabel}[1]{\label{coro:#1}}
\nc{\coref}[1]{Corollary~\ref{coro:#1}}
\nc{\exlabel}[1]{\label{exam:#1}}
\nc{\exref}[1]{Example~\ref{exam:#1}}
\nc{\delabel}[1]{\label{defi:#1}}
\nc{\deref}[1]{Definition~\ref{defi:#1}}
\nc{\eqlabel}[1]{\label{equation:#1}}
\nc{\eqref}[1]{(\ref{equation:#1})}
\nc{\csm}{\mbox{$\triangleright\!\!\!<$}}
\nc{\smc}{\mbox{$>\!\!\!\triangleleft$}}
\nc{\trr}{\triangleright}
\providecommand{\operatorname}[1]{\mathrm{#1}\,}
\nc{\Hom}{\operatorname{Hom}}
\nc{\Mor}{\operatorname{Mor}}
\nc{\Aut}{\operatorname{Aut}}
\nc{\Ann}{\operatorname{Ann}}
\nc{\Ker}{\operatorname{Ker}}
\nc{\Trace}{\operatorname{Trace}}
\nc{\Char}{\operatorname{Char}}
\nc{\Mod}{\operatorname{Mod}}
\nc{\End}{\operatorname{End}}
\nc{\Spec}{\operatorname{Spec}}
\nc{\Span}{\operatorname{Span}}
\nc{\sgn}{\operatorname{sgn}}
\nc{\Id}{\operatorname{Id}}
\nc{\Com}{\operatorname{Com}}

\def\Box{\mbox{$\sqcap\!\!\!\!\sqcup$}}

\nc{\dht}{\mbox{$\rightharpoonup\hspace{-2ex}\rightharpoonup$}}
\nc{\dhtb}{\mbox{$\leftharpoonup\hspace{-2ex}\leftharpoonup$}}
\nc{\nd}{\mbox{$\not|$}} %not divide sign
\providecommand{\text}[1]{\mbox{{\textrm #1}}}

\nc{\nci}{\mbox{$\not\subseteq$}}
\nc{\scontainin}{\mbox{$\mbox{}\subseteq\hspace{-1.5ex}\raisebox{-.5ex}{$_\prime
$}\hspace*{1.5ex}$}}

\def\ot{\otimes}

\def\doublerightleft#1#2{{\lower.2ex\vbox{
\hbox{${\smash{\mathop{\longrightarrow}\limits^{#1}}}$}\vspace*{-4mm}
\hbox{${\smash{\mathop{\longleftarrow}\limits_{#2}}}$}}}}

\newfont{\bbb}{msbm10 scaled\magstep1}  % Blackboardbold for 12pt article
\newfont{\bbbsub}{msbm10}                % Blackboardbold for subscripts for                                     % 12pt article
\newfont{\msam}{msam10 scaled\magstep1}

\begin{document}
\title{A class of non-symmetric solutions for the integrability condition
of the Knizhnik-Zamolodchikov equation: a Hopf algebra approach}
\author{G. Militaru
\\University of Bucharest
\\Faculty of Mathematics\\Str. Academiei 14
\\RO-70109 Bucharest 1, Romania
\\e-mail: gmilit@al.math.unibuc.ro}
\date{}
\maketitle
\begin{abstract}
\noindent
Let $M$ be a vector space over a field $k$ and $R\in \End_k(M\ot M)$.
This paper studies what shall be called
the Long equation: that is, the system of nonlinear equations
$R^{12}R^{13}=R^{13}R^{12}$ and $R^{12}R^{23}=R^{23}R^{12}$ in
$\End_k(M\ot M\ot M)$. Any symmetric solution of this system supplies
us a solution of the integrability condition of the
Knizhnik-Zamolodchikov equation: $[R^{12}, R^{13}+R^{23}]=0$
(\cite{K} or \cite{SS}).
We shall approach this equation by introducing a new class of bialgebras,
which we call Long bialgebras: these are pairs
$(H,\sigma)$, where $H$ is a bialgebra and $\sigma :H\ot H\to k$ is a
$k$-bilinear map satisfying certain properties. The main theorem of
this paper is a FRT type theorem: if $M$ is finite dimensional, any
solution $R$ of the Long equation has the form $R=R_{\sigma}$, where
$M$ has a structure of a right comodule over a Long bialgebra
$(L(R), \sigma)$, and $R_{\sigma}$ is the special map
$R_{\sigma}(m\ot n)=\sum \sigma(m_{<1>}\ot n_{<1>})m_{<0>}\ot n_{<0>}$.
\end{abstract}

\section{Introduction}
Let $M$ be a finite dimensional vector space over a field $k$ and
$R\in \End_k(M\ot M)$. In this paper we shall describe all solutions
of what we have called the Long equation, namely the system of nonlinear
equations
\begin{equation}\label{LO}
\left\lbrace
\begin{array}{ccc}
R^{12}R^{13}=R^{13}R^{12}\\
R^{12}R^{23}=R^{23}R^{12}
\end{array}
\right.
\end{equation}
in $\End_k(M\ot M\ot M)$. Our approach is similar to the one used by
Faddeev, Reshetikhin and Takhtadjian in relating the
solutions of the quantum Yang-Baxter equation to comodules over
co-quasitriangular bialgebras (see \cite{FRT}).

We shall introduce a new class of bialgebras which we shall
call Long bialgebras. They are pairs $(H,\sigma)$, where
$H$ is a bialgebra and $\sigma:H\ot H\to k$ is a linear map satisfying
the conditions $(L1)-(L5)$ from definition \ref{LON}. The conditions
$(L2)-(L5)$ are identical with the conditions $(B2)-(B5)$ from the
definition of co-quasitriangular bialgebras. What differentiates
Long bialgebras from the co-quasitriangular bialgebras is the condition
$(L1)$ versus $(B1)$. This new class of bialgebras will play a
fundamental role in solving the Long equation. More precisely,
if $(M,\rho)$ is a right comodule over a Long bialgebra $(H,\sigma)$,
then the special map
$$
R_{\sigma}:M\ot M\to M\ot M,\quad
R_{\sigma}(m\ot n)=\sum \sigma(m_{<1>}\ot n_{<1>})m_{<0>}\ot n_{<0>}
$$
is a solution of the Long equation. Conversely, the main theorem of this
paper is a FRT type theorem: if $M$ is a finite dimensional vector space
and $R$ is a solution of the Long equation, then there exists a
Long bialgebra $(L(R), \sigma)$ such that $M$ has a structure of
right $L(R)$-comodule and  $R=R_{\sigma}$.
Let us now look at the Long equation from a different
angle, having in mind Radford's version of the FRT theorem (\cite{R}):
in the finite dimensional case, any solution $R$ of the quantum Yang-Baxter
equation has the form $R=R_{(M,\cdot,\rho)}$, where
$(M,\cdot,\rho)\in {}_{A(R)}{\cal YD}^{A(R)}$, the category of
Yetter-Drinfel'd modules. This is obtained immediately, keeping in mind
the fact that a comodule $(M,\rho)$ over a co-quasitriangular bialgebra
(like $(A(R),\sigma )$) has a structure of Yetter-Drinfel'd module via
$$
h\cdot m=\sum \sigma(m_{<1>}\ot h)m_{<0>}
$$
for all $h\in A(R)$, $m\in M$ and $R_{(M,\cdot,\rho)}=R_{\sigma}$.
A similar phenomenon happens to the comodules over a Long bialgebra
$(H,\sigma)$: they become objects in the category
${}_H{\cal L}^H$ of $H$-dimodules (see Proposition 3.6). This category has
been introduced by Long (\cite{L}) for the case of a commutative and
cocommutative $H$ and studied in connection with the Brauer group of an
$H$-dimodule algebra. For this reason, we called system (\ref{LO})
the Long equation and the new class of bialgebras, Long bialgebras.

Finally, we shall point out a connection between the symmetric solutions
of the Long equation (e.g. the map $R^{\phi}$ from proposition \ref{vre})
and solutions of the integrability condition of the
Knizhnik-Zamolodchikov equation. Suppose that $R$ is a solution of the
Long equation; then $R$ satisfies the equation
$[R^{12}, R^{13}+R^{23}]=0$,
which is called the integrability condition of the Knizhnik-Zamolodchikov
equation (see \cite{SS}). Moreover, we assume that $R$ is symmetric,
that is $R^{12}=R^{21}$. $W$ is a solution of the Knizhnik-Zamolodchikov
equation if it is a solution of the following system of differential
equations:
\begin{equation}\label{doi}
\frac{\partial W}{\partial z^i}=h\:\sum_{i\neq j}
(\frac{R^{ij}}{z^i-z^j})W.
\end{equation}
This is the equation for a function $W(z)$ taking values in $M^{\ot n}$,
representing a covariant constant section of the trivial bundle
$Y_n \times M^{\ot n}\to Y_n$ with flat connection
$\sum_{i\neq j}R^{ij} \frac{dz^i}{z^i-z^j}$ (here $k$ is the complex field
$\mathbf{C}$, $h$ is a complex parameter and
$Y_n=\mathbf{C}^n/\mbox{multidiagonal}$).
In other words, $R$ is describing a family of flat connections on bundles
with fiber $M^{\ot n}$. For further details we refer to \cite{K} and
\cite{SS}.

\section{Preliminaries}
Throughout this paper, $k$ will be a field.
All vector spaces, algebras, coalgebras and bialgebras considered
are over $k$. $\ot$ and $\Hom$ will mean $\ot_k$ and $\Hom_k$.
For a coalgebra $C$, we will use Sweedler's $\Sigma$-notation, that is,
$\Delta(c)=\sum c_{(1)}\ot c_{(2)},~(I\ot\Delta)\Delta(c)=
\sum c_{(1)}\ot c_{(2)}\ot c_{(3)}$, etc. We will
also use  Sweedler's notation for right $C$-comodules:
$\rho_M(m)=\sum m_{<0>}\otimes m_{<1>}$, for any $m\in M$ if
$(M,\rho_M)$ is a right $ C$-comodule. ${\cal M}^C$ will be the
category of right $C$-comodules and $C$-colinear maps and
${}_A{\cal M}$ will be the category of left $A$-modules and
$A$-linear maps, if $A$ is a $k$-algebra. An important role in the present
paper will be played by ${\cal M}^n(k)$, the comatrix coalgebra of
order $n$, i.e. ${\cal M}^n(k)$ is the $n^2$-dimensional vector space
with $\{c_{ij}\mid i,j=1,\cdots,n \}$ a $k$-basis such that
\begin{equation}\label{com1}
\Delta(c_{jk})=\sum_{u=1}^{n}c_{ju}\ot c_{uk},
\quad \varepsilon (c_{jk})=\delta_{jk}
\end{equation}
for all $j,k=1,\cdots, n$. We view $T({\cal M}^n(k))$ with the unique
bialgebra structure which can be defined on the tensor algebra
$T({\cal M}^n(k))$, which extends the comultiplication $\Delta$ and the
counity $\varepsilon$ of ${\cal M}^n(k)$.

Let $H$ be a bialgebra. An {\sl $H$-dimodule } over $H$ is a
triple $(M,\cdot ,\rho)$, where $(M,\cdot)$ is a left $H$-module,
$(M,\rho)$ is a right $H$-comodule such that the following
compatibility condition holds:
\begin{equation}\label{C}
\rho(h\cdot m)=\sum h\cdot m_{<0>}\ot m_{<1>}
\end{equation}
for all $h\in H$ and $m\in M$.
The category of $H$-dimodules over $H$ and $H$-linear $H$-colinear
maps will be denoted by ${}_H{\cal L}^H$. This category was introduced
for a commutative and cocommutative $H$ by Long in \cite{L}.

For a vector space $M$, $\tau :M\otimes M\to M\otimes M$
will denote the flip map, that is $\tau (m\otimes n)=n\otimes m$
for all $m$, $n \in M$. If $R:M\ot M\to M\ot M$ is a linear map,
we denote by $R^{12}$, $R^{13}$, $R^{23}$ the maps of $\End_k(M\ot M\ot M)$
given by
$$
R^{12}=R\ot I, \quad R^{23}=I\ot R,\quad
R^{13}=(I\ot \tau)(R\ot I)(I\ot \tau).
$$

\section{The Long equation}
We shall start with the following:

\begin{definition}
Let $M$ be a vector space and $R\in \End_k(M\ot M)$. We shall say that
$R$ is a solution for the Long equation if
\begin{equation}\label{ICKZ}
\left\lbrace
\begin{array}{ccc}
R^{12}R^{13}=R^{13}R^{12}\\
R^{12}R^{23}=R^{23}R^{12}
\end{array}
\right.
\end{equation}
holds in $\End_k(M\ot M\ot M)$.
\end{definition}

\begin{remarks}
1. If $R$ is a solution of the Long equation, then $R$ satisfies the
equation
$$[R^{12}, R^{13}+R^{23}]=0.$$
Hence, any symmetric solution $R$ (i.e. $R^{12}=R^{21}$) of the
Long equation is a solution for the integrability condition
of the Knizhnik-Zamolodchikov equation (see \cite{K} or \cite{SS}).

2. Let $M$ be a finite dimensional vector space and $\{m_1,\cdots,m_n \}$
a basis of $M$. Let  $R\in \End_k(M\ot M)$ given by
$$
R(m_v\ot m_u)=\sum_{i,j}x_{uv}^{ji}m_i\ot m_j,
$$
for all $u$, $v=1,\cdots ,n$, where $(x_{uv}^{ji})_{i,j,u,v}$
is a family of scalars of $k$. Then $R$ is a solution of the
Long equation if and only if the following two equations hold:
\begin{equation}\label{cinci}
\sum_{v}x_{kv}^{ji}x_{ql}^{pv}=
\sum_{\alpha}x_{kl}^{j\alpha}x_{q\alpha }^{pi}
\end{equation}
\begin{equation}\label{sase}
\sum_{v}x_{kv}^{ji}x_{lq}^{vp}=
\sum_{\alpha}x_{kl}^{j\alpha}x_{\alpha q}^{ip}
\end{equation}
for all $i$, $j$, $k$, $l$, $p$, $q=1,\cdots, n$.

3. The second equation of our system, namely
$R^{12}R^{23}=R^{23}R^{12}$ (called the ${\cal D}$-equation),
which is obtained from the quantum Yang-Baxter equation
$R^{12}R^{13}R^{23}=R^{23}R^{13}R^{12}$
by deleting the middle term from both sides, was studied in \cite{M3}.
For a bialgebra $H$, let ${}_H{\cal L}^H$ be the category of $H$-dimodules.
In \cite{M3} we proved that in the finite dimensional case the maps
$$
R_{(M,\cdot,\rho)}:M\ot M\to M\ot M,\quad
R_{(M,\cdot,\rho)}(m\ot n):=\sum n_{<1>}\cdot m\ot n_{<0>},
$$
where $(M,\cdot,\rho)$ is a $D(R)$-dimodule for some bialgebra $D(R)$,
describle all the solutions of the ${\cal D}$-equation.
\end{remarks}

\begin{examples}
1. Suppose that $R\in \End_k(M\ot M)$ is bijective. Then, $R$ is a
solution of the Long equation if and only if $R^{-1}$ is.

2. Let $(m_i)_{i\in I}$ be a basis of $M$ and $(a_{ij})_{i,j\in I}$ be a
family of scalars of $k$. Then, $R:M\ot M\to M\ot M$,
$R(m_i\ot m_j)=a_{ij}m_i\ot m_j$, for all $i$, $j\in I$, is a solution of
the Long equation. In particular, the identity map $Id_{M\ot M}$ is a
solution of the Long equation.

3. Let $M$ be a finite dimensional vector space and $u$ an automorphism
of $M$. If $R$ is a solution of the Long equation, then
${}^{u}R:= (u\ot u)R(u\ot u)^{-1}$ is also a solution of the
Long equation.

4. Let $A$ be a $k$-algebra, $(a_i)_{i=1,\cdots,n}$ be a family of elements
of $A$ and $(M,\cdot)$ a left $A$-module. Then
$$
R:M\ot M\to M\ot M,\quad R(l\ot m):=\sum_{i=1}^{n}l\ot a_i\cdot m
$$
for all $l$, $m\in M$, is a solution of the Long equation.

5. Let $A$ be a $k$-algebra and $R=\sum R^1\ot R^2 \in A\ot A$ such that
\begin{equation}\label{cen}
\sum R^1a\ot R^2=\sum aR^1\ot R^2
\end{equation}
for all $a\in A$.
If $(M,\cdot)$ is a left $A$-module, then the homothety
$$
{\cal R}:M\ot M\to M\ot M,\quad {\cal R}(m\ot n)=R\cdot (m\ot n)=
\sum R^1\cdot m\ot R^2\cdot n,
$$
is a solution of the Long equation.

6. Let $f$, $g\in \End_k(M)$ such that $fg=gf$. Then, $R:=f\ot g$
is a solution of the Long equation.

In particular, let $M$ be a two-dimensional vector space with
$\{m_1, m_2 \}$ a basis. Let $f$, $g\in \End_k(M)$ such that, with
respect to the given basis, they are:
$$
f=
\left(
\begin{array}{cc}
a&1\\
0&a
\end{array}
\right), \quad
g=
\left(
\begin{array}{cc}
b&c\\
0&b
\end{array}
\right)
$$
where $a$, $b$, $c$ are scalars of $k$. Then, $R=f\ot g$, with respect
to the ordered basis $\{m_1\ot m_1, m_1\ot m_2, m_2\ot m_1, m_2\ot m_2 \}$
of $M\ot M$, is given by
\begin{equation}\label{100}
R=
\left(
\begin{array}{cccc}
ab&ac&b&c\\
0&ab&0&b\\
0&0&ab&ac\\
0&0&0&ab
\end{array}
\right)
\end{equation}
and $R$ is a solution of the Long equation.

7. Let $G$ be an abelian group and $(M,\cdot)$ be a $k[G]$-module.
Suppose that there exists
$\{M_{\sigma }\mid \sigma \in G \}$ a family of $k[G]$-submodules
of $M$ such that $M=\oplus_{\sigma \in G}M_{\sigma}$.
If $m\in M$, then $m$ is a finite sum of homogenous elements
$m=\sum m_{\sigma}$. The map
\begin{equation}\label{gr}
R:M\ot M\to M\ot M, \quad R(n\ot m)=
\sum_{\sigma}\sigma\cdot n\ot m_{\sigma},\; \forall n, m\in M
\end{equation}
is a solution of the Long equation.
\end{examples}

In \cite{IS}, for an $n$-dimensional vector space $M$, a certain operator
$R_{\phi}:M\ot M\to M\ot M$ is associated to any function
$\phi:\{1,\cdots,n \}\to \{1,\cdots,n \}$, with $\phi^{2}=\phi$.
This operator  $R_{\phi}$ is a solution for what we called in
\cite{M1} the Hopf equation
$$
R^{23}R^{13}R^{12}=R^{12}R^{23}
$$
We shall modify the operator $R_{\phi}$ in order to make it a solution
for the Long equation.

\begin{proposition}\label{vre}
Let $M$ be an $n$-dimensional vector space with $\{m_1,\cdots, m_n\}$
a basis and $\phi:\{1,\cdots,n \}\to \{1,\cdots,n \}$ a function
with $\phi^{2}=\phi$. Let
$$
R^{\phi}:M\ot M\to M\ot M,\quad
R^{\phi}(m_i\ot m_j)=\delta_{ij}\delta_{i\in Im(\phi)}
\sum_{a,b\in \phi^{-1}(i)}m_a\ot m_b
$$
for all $i$, $j=1,\cdots, n$. Then
$R:=R^{\phi}$ is a symmetric solution for the equation
$$
R^{12}R^{13}=R^{13}R^{12}=R^{12}R^{23}=R^{23}R^{12}.
$$
In particular, $R^{\phi}$ is a solution of the Long equation.
\end{proposition}

\begin{proof} For $i$, $j$, $k=1,\cdots,n$ we have:
\begin{eqnarray*}
R^{12}R^{13}(m_i\ot m_j\ot m_k)&=&
R^{12}\Bigl (\delta_{ik}\delta_{i\in Im(\phi)}
\sum_{a,b\in \phi^{-1}(i)}m_a\ot m_j\ot m_b \Bigl )\\
&=&\delta_{ik}\delta_{i\in Im(\phi)}\sum_{a,b\in \phi^{-1}(i)}
\delta_{aj}\delta_{a\in Im(\phi)}\sum_{c,d\in \phi^{-1}(a)}
m_c\ot m_d\ot m_b\\
&=&\delta_{ik}\delta_{i\in Im(\phi)}
\delta_{j\in \phi^{-1}(i)}\delta_{j\in Im(\phi)}
\sum_{b\in \phi^{-1}(i), c,d\in \phi^{-1}(j)}
m_c\ot m_d\ot m_b\\
\text{(\mbox{using} $\phi^2=\phi$ )}
&=&\delta_{ik}\delta_{ij}\delta_{i\in \phi^{-1}(i)}
\sum_{b,c,d\in \phi^{-1}(i)}m_c\ot m_d\ot m_b
\end{eqnarray*}
and
\begin{eqnarray*}
R^{13}R^{12}(m_i\ot m_j\ot m_k)&=&
R^{13}\Bigl (\delta_{ij}\delta_{i\in Im(\phi)}
\sum_{a,d\in \phi^{-1}(i)}m_a\ot m_d\ot m_k \Bigl )\\
&=&\delta_{ij}\delta_{i\in Im(\phi)}\sum_{a,d\in \phi^{-1}(i)}
\delta_{ak}\delta_{a\in Im(\phi)}\sum_{c,b\in \phi^{-1}(a)}
m_c\ot m_d\ot m_b\\
&=&\delta_{ij}\delta_{i\in Im(\phi)}
\delta_{k\in \phi^{-1}(i)}\delta_{k\in Im(\phi)}
\sum_{d\in \phi^{-1}(i), b,c\in \phi^{-1}(k)}
m_c\ot m_d\ot m_b\\
\text{(\mbox{using} $\phi^2=\phi$ )}
&=&\delta_{ij}\delta_{ik}\delta_{i\in \phi^{-1}(i)}
\sum_{b,c,d\in \phi^{-1}(i)}m_c\ot m_d\ot m_b
\end{eqnarray*}
i.e. $R^{12}R^{13}=R^{13}R^{12}$. The proof of the second identity
($R^{12}R^{23}=R^{23}R^{12}$) is left to the reader.

On the other hand the family of scalars $(x_{uv}^{ji})$ which define
$R^{\phi}$ are given by
$x_{uv}^{ji}=\delta_{uv}\delta_{\phi(i)v}\delta_{\phi(j)v}$
for all $i$, $j$, $u$, $v=1,\cdots ,n$. Hence, $x_{uv}^{ji}=x_{vu}^{ji}$,
i.e. $R$ is symmetric.
\end{proof}

\section{Long bialgebras }
Recall from \cite{M3} the following definition

\begin{definition}
Let $C$ be a coalgebra. A $k$-bilinear map
$\sigma :C\ot C\to k$ is called a strong ${\cal D}$-map if
\begin{equation}\label{spe}
\sum \sigma(c_{(1)}\ot d)c_{(2)}=
\sum \sigma(c_{(2)}\ot d)c_{(1)}
\end{equation}
for all $c$, $d\in C$.
\end{definition}

\begin{proposition}\label{muc}
Let $C$ be a coalgebra and $\sigma :C\ot C\to k$ a strong ${\cal D}$-map.
Let $(M,\rho)$ be a right $C$-comodule. Then, the special map
$$
R_{\sigma}:M\ot M\to M\ot M, \quad
R_{\sigma}(m\ot n)=
\sum \sigma(m_{<1>}\ot n_{<1>})m_{<0>}\ot n_{<0>}
$$
is a solution of the Long equation.
\end{proposition}

\begin{proof} Let $R=R_{\sigma}$. Then, the fact that
$R^{12}R^{23}=R^{23}R^{12}$ is Proposition 4.3 of \cite{M3}.
We show that $R^{12}R^{13}=R^{13}R^{12}$. For $l$, $m$, $n\in M$ we have:
\begin{eqnarray*}
R^{12}R^{13}(l\ot m\ot n)&=&
R^{12}\Bigl(\sum \sigma(l_{<1>}\ot n_{<1>})
l_{<0>}\ot m\ot n_{<0>}\Bigl)\\
&=&\sum \sigma(l_{<2>}\ot n_{<1>})
\sigma(l_{<1>}\ot m_{<1>})l_{<0>}\ot m_{<0>}\ot n_{<0>}
\end{eqnarray*}
and
\begin{eqnarray*}
R^{13}R^{12}(l\ot m\ot n)&=&
R^{13}\Bigl(\sum \sigma(l_{<1>}\ot m_{<1>})
l_{<0>}\ot m_{<0>}\ot n\Bigl)\\
&=&\sum \sigma(l_{<2>}\ot m_{<1>})
\sigma(l_{<1>}\ot n_{<1>})l_{<0>}\ot m_{<0>}\ot n_{<0>}\\
&=&\sum \sigma \Bigl( \underline{\sigma(l_{<1>(1)}\ot n_{<1>})l_{<1>(2)}}
\ot m_{<1>}\Bigl)l_{<0>}\ot m_{<0>}\ot n_{<0>}\\
\text{(\mbox{using (\ref{spe})})}
&=&\sum \sigma \Bigl( \underline{\sigma(l_{<1>(2)}\ot n_{<1>})l_{<1>(1)}}
\ot m_{<1>}\Bigl)l_{<0>}\ot m_{<0>}\ot n_{<0>}\\
&=&\sum \sigma(l_{<2>}\ot n_{<1>})
\sigma(l_{<1>}\ot m_{<1>})l_{<0>}\ot m_{<0>}\ot n_{<0>}
\end{eqnarray*}
i.e. $R^{12}R^{13}=R^{13}R^{12}$. Hence, $R_{\sigma}$ is a solution of
the Long equation.
\end{proof}

Now, we shall introduce a new class of bialgebras which play for the
Long equation the same role as the co-quasitriangular (or braided)
bialgebras do for the quantum Yang-Baxter equation.

\begin{definition}\label{LON}
A Long bialgebra is a pair $(H,\sigma)$, where $H$ is a bialgebra and
$\sigma :H\ot H\to k$ is a $k$-linear map such that the following
conditions are fulfilled:

$(L1)\quad \sum \sigma(x_{(1)}\ot y)x_{(2)}=
\sum \sigma(x_{(2)}\ot y)x_{(1)}$

$(L2)\quad \sigma (x\ot 1)=\varepsilon(x)$

$(L3) \quad \sigma(x\ot yz)=
\sum \sigma (x_{(1)}\ot y)\sigma (x_{(2)}\ot z)$

$(L4)\quad \sigma (1\ot x)=\varepsilon(x)$

$(L5) \quad \sigma(xy\ot z)=
\sum \sigma (y\ot z_{(1)})\sigma (x\ot z_{(2)})$

for all $x$, $y$, $z\in H$.
\end{definition}

\begin{remark}
We recall that (see, for instance, \cite{K}) a pair $(H,\sigma)$ is a
co-quasitriangular (or braided) bialgebra if $\sigma$ satisfies the
conditions $(B1)$-$(B5)$. The conditions $(B2)$-$(B5)$ are identical with the
conditions $(L2)$-$(L5)$. The first condition is

$(B1)\quad \sum \sigma(x_{(1)}\ot y_{(1)})y_{(2)}x_{(2)}=
\sum \sigma(x_{(2)}\ot y_{(2)})x_{(1)}y_{(1)}$

for all $x$, $y\in H$. If $(H,\sigma)$ is a co-quasitriangular
bialgebra and $(M,\rho)$ is a right $H$-comodule, then the special
map
$$
R_{\sigma}:M\ot M\to M\ot M, \quad
R_{\sigma}(m\ot n)=\sum \sigma(m_{<1>}\ot n_{<1>})m_{<0>}\ot n_{<0>}
$$
is a solution for the quantum Yang-Baxter equation
$R^{12}R^{13}R^{23}=R^{23}R^{13}R^{12}$.
Conversely, if $M$ is a finite dimensional vector space and $R$ is a
solution of the quantum Yang-Baxter equation, then there exists a
bialgebra $A(R)$ and a unique $k$-bilinear map $\sigma :A(R)\ot A(R)\to k$
such that $(A(R), \sigma)$ is co-quasitriangular,
$M\in {\cal M}^{A(R)}$ and $R=R_{\sigma}$ (see \cite{FRT} or \cite{K}).
\end{remark}

We shall now present a few examples of Long bialgebras. More examples will
be given after the main result.

\begin{examples}
1. If $H$ is cocommutative, then $(L1)$ holds for any $k$-linear map
$\sigma :H\ot H\to k$. In particular, if $G$ is a group and
$\sigma :G\times G\to k$ is a bicharacter on $G$, then $(k[G], \sigma)$ is
a Long bialgebra.

2. Let $H=k<x,y>$ be the free algebra generated by $x$ and $y$, with the
bialgebra structure given by
$$
\Delta(x)=x\ot x,\quad \Delta(y)=y\ot 1+x\ot y,\quad
\varepsilon(x)=1,\quad \varepsilon(y)=0.
$$
Let $\sigma :H\ot H\to k$ be a $k$-bilinear map. Then, it is easy to
see that $(H,\sigma)$ is a Long bialgebra if and only if
$$
\begin{array}{lll}
\sigma(x\ot 1)=1, &\sigma(x\ot x)=1, & \sigma(x\ot y)=0, \\
\sigma(y\ot 1)=0, &\sigma(y\ot x)=0, & \sigma(y\ot y)=0.
\end{array}
$$
3. Let $H=k<x,y,z>$ be the free algebra generated by $x$, $y$, $z$,
with the bialgebra structure given by
$$
\Delta(x)=x\ot x,\quad \Delta(y)=y\ot y,\quad \Delta(z)=x\ot z+z\ot y
$$
$$
\varepsilon(x)=\varepsilon(y)=1,\quad \varepsilon(z)=0.
$$
Let $\sigma :H\ot H\to k$ be a $k$-bilinear map. Then
$(H,\sigma)$ is a Long bialgebra if and only if there exist $a$, $b$,
$c$ scalars of $k$ such that
$$
\sigma(x\ot x)=\sigma(y\ot x)=a,\quad \sigma(x\ot y)=\sigma(y\ot y)=b,\quad
\sigma(x\ot z)=\sigma(y\ot z)=c,
$$
$$
\sigma(z\ot x)=\sigma(z\ot y)=\sigma(z\ot z)=0.
$$
4. Let $H_4$ be Sweedler's $4$-dimensional Hopf algebra, i.e.
$$
H_4=k<x,y,z\mid \; x^2=1,\; y^2=0,\; xy=z,\; xz=-zx=y>
$$
with
$$
\Delta(x)=x\ot x,\quad \Delta(y)=y\ot x+1\ot y,\quad
\Delta(z)=x\ot z+z\ot 1
$$
$$
\varepsilon(x)=1,\quad \varepsilon(y)=\varepsilon(z)=0.
$$
It is well known that $H_4$ is a co-quasitriangular bialgebra
(see, for example, \cite{Doi}). We shall prove that
there exists no $\sigma :H_4\ot H_4\to k$ such that
$(H_4,\sigma)$ is a Long bialgebra.

Suppose that there exists $\sigma :H_4\ot H_4\to k$ such that
$(H_4,\sigma)$ is a Long bialgebra. It follows from $(L1)$ that
$$
\sigma(y\ot h)x+\sigma(1\ot h)y=\sigma(x\ot h)y+\sigma(y\ot h)1_H
$$
for all $h\in H$. As $\{1, x, y, z\}$ is a basis of $H_4$, we get that
$\sigma(x\ot h)=\sigma(1\ot h)$ and $\sigma(y\ot h)=0$
for all $h\in H$. In particular,
$\sigma(x\ot x)=1$ and $\sigma(x\ot z)=\sigma(x\ot y)=0$.
It follows that
$$
0=\sigma(x\ot y)=\sigma(x\ot xz)=\sigma(x\ot x)+\sigma(x\ot z)
$$
hence $\sigma(x\ot z)=-1$, contradiction.
\end{examples}

\begin{proposition}\label{gen}
Let $H$ be a bialgebra and $\sigma :H\ot H\to k$ a $k$-bilinear map
which satisfies $(L3)$ and $(L5)$. Suppose that $(L1)$ holds for
a system of generators of $H$ as an algebra. Then $(L1)$ holds for all
elements of $H$.
\end{proposition}

\begin{proof}
Let $x$, $y$, $z\in H$ be three elements among the generators of $H$.
It is enough to prove that $(L1)$ holds for $(x,yz)$ and $(xy,z)$.
We have:
\begin{eqnarray*}
\sum \sigma (x_{(1)}\ot yz)x_{(2)}&=&
\sum \sigma (x_{(1)(1)}\ot y) \sigma (x_{(1)(2)}\ot z)x_{(2)}\\
&=&\sum \sigma (x_{(1)}\ot y)
\underline{\sigma (x_{(2)}\ot z)x_{(3)}}\\
\text{(\mbox{using (L1)} )}
&=&\sum\sigma (x_{(1)}\ot y)\sigma (x_{(3)}\ot z)x_{(2)}\\
&=&\sum\underline{\sigma (x_{(1)}\ot y)x_{(2)}}\sigma (x_{(3)}\ot z)\\
\text{(\mbox{using (L1)} )}
&=&\sum \sigma (x_{(2)}\ot y)\sigma (x_{(3)}\ot z)x_{(1)}\\
\text{(\mbox{using (L3)} )}
&=&\sum \sigma (x_{(2)}\ot yz)x_{(1)}
\end{eqnarray*}
i.e. $(L1)$ holds for $(x,yz)$. On the other hand
\begin{eqnarray*}
\sum \sigma ((xy)_{(1)}\ot z)(xy)_{(2)}&=&
\sum \sigma (x_{(1)}y_{(1)}\ot z) x_{(2)}y_{(2)}\\
&=&\sum \sigma (y_{(1)}\ot z_{(1)})
\underline{\sigma (x_{(1)}\ot z_{(2)})x_{(2)}}y_{(2)}\\
\text{(\mbox{using (L1)} )}
&=&\sum\sigma (x_{(2)}\ot z_{(2)})x_{(1)}
\underline{\sigma (y_{(1)}\ot z_{(1)})y_{(2)}}\\
\text{(\mbox{using (L1)} )}
&=&\sum\sigma (x_{(2)}\ot z_{(2)})x_{(1)}
\sigma (y_{(2)}\ot z_{(1)})y_{(1)}\\
\text{(\mbox{using (L5)} )}
&=&\sum \sigma (x_{(2)}y_{(2)}\ot z)x_{(1)}y_{(1)}\\
&=&\sum \sigma ((xy)_{(2)}\ot z)(xy)_{(1)}
\end{eqnarray*}
and the proof is complete now.
\end{proof}

We recall that if $M$ is a right $H$-comodule over a co-quasitriangular
bialgebra $(H,\sigma)$, then $M$ can be viwed as an object of
${}_H{\cal YD}^H$, the category of Yetter-Drinfel'd modules,
where the structure of left $H$-module on $M$ is induced by $\sigma$:
$$
h\cdot m=\sum \sigma(m_{<1>}\ot h)m_{<0>}
$$
for all $h\in H$, $m\in M$. In the next proposition
we shall prove that any right $H$-comodule $M$ over a
Long bialgebra $(H,\sigma)$ can be viewed as an $H$-dimodule.

\begin{proposition}\label{bic}
Let $(H,\sigma)$ be a Long bialgebra and $(M,\rho)$ be a right
$H$-comodule. Then, the left action of $H$ on $M$ given by
$$
h\cdot m=\sum \sigma (m_{<1>}\ot h)m_{<0>}
$$
for all $h\in H$, $m\in M$, makes $(M,\cdot,\rho)$ an $H$-dimodule.
\end{proposition}

\begin{proof}
The fact that $(M,\cdot)$ is a left $H$-module follows from $(L2)$ and
$(L3)$. We shall prove that the compatibility condition (\ref{C}) holds.
For $h\in H$ and $m\in M$ we have
\begin{eqnarray*}
\rho(h\cdot m)&=&\rho\Bigl(\sum \sigma (m_{<1>}\ot h)m_{<0>}\Bigl)\\
&=&\sum \sigma (m_{<1>}\ot h)m_{<0><0>}\ot m_{<0><1>}\\
&=&\sum m_{<0>}\ot \sigma(m_{<1>(2)}\ot h)m_{<1>(1)}\\
\text{(\mbox{using (L1)} )}
&=&\sum m_{<0>}\ot \sigma(m_{<1>(1)}\ot h)m_{<1>(2)}\\
&=&\sum \sigma (m_{<1>}\ot h)m_{<0>}\ot m_{<2>}\\
&=&\sum h\cdot m_{<0>}\ot m_{<1>}
\end{eqnarray*}
i.e. $(M,\cdot,\rho)$ is an $H$-dimodule.
\end{proof}

From Proposition \ref{muc} we obtain that if $(H,\sigma)$ is a
Long bialgebra and $M$ is a right $H$-comodule, then the natural map
$R_{\sigma}$ is a solution of the Long equation. Now, we shall prove
the main result of this paper: in the finite dimensional case any solution
of the Long equation arises in this way.

\begin{theorem}\label{boci}
Let $M$ be a finite dimensional vector space and $R\in \End_k(M\ot M)$
be a solution for the Long equation.
Then there exists a bialgebra $L(R)$ and a unique $k$-bilinear map
$$\sigma:L(R)\ot L(R)\to k$$
such that $(L(R), \sigma)$ is a Long bialgebra, $M$ is a right
$L(R)$-comodule and  $R=R_{\sigma}$.

Furthermore, if $R$ is bijective, then $\sigma$ is invertible in
convolution.
\end{theorem}

\begin{proof}
Suppose that dim$_k(M)=n$ and let  $\{m_1,\cdots ,m_n \}$
be a basis for $M$. Let $(x_{uv}^{ji})_{i,j,u,v}$ be a family
of scalars of $k$ such that
\begin{equation}
R(m_v\ot m_u)=\sum_{i,j}x_{uv}^{ji}m_i\ot m_j
\end{equation}
for all $u$, $v=1,\cdots ,n$.

Let ${\cal M}^n(k)$, be the comatrix coalgebra
of order $n$ with $\{c_{ij}\mid i, j=1, \cdots ,n \}$ a $k$-basis and
$(T({\cal M}^n(k)),\Delta,\varepsilon, \mu, 1)$ the unique
bialgebra structure on the tensor algebra $(T({\cal M}^n(k)), \mu, 1)$ whose
comultiplication $\Delta$ and counit $\varepsilon$ extend the
comultiplication and the counit from ${\cal M}^n(k)$. Then $M$ has a unique
right $T({\cal M}^n(k))$-comodule structure
$\rho :M\to M\ot T({\cal M}^n(k))$ such that
\begin{equation}\label{ro}
\rho(m_l)=\sum_{v=1}^{n}m_v\ot c_{vl}
\end{equation}
for all $l=1,\cdots, n$. Let $o(i,j,k,l)$ be the obstruction
\begin{equation}\label{obst}
o(i,j,k,l):=\sum_{v}x_{kv}^{ji}c_{vl} -
\sum_{\alpha}x_{kl}^{j\alpha}c_{i\alpha}
\end{equation}
for all $i$, $j$, $k$, $l=1,\cdots, n$. Let $I(R)$ be the two-sided ideal
of $T({\cal M}^n(k))$ generated by all $o(i,j,k,l)$. Then, $I(R)$ is also a
coideal of $T({\cal M}^n(k))$, i.e. it is a biideal. This follows from the
formulas
$$
\Delta(o(i,j,k,l))=\sum_{u}\Bigl( o(i,j,k,u)\ot c_{ul}+
c_{iu}\ot o(u,j,k,l)\Bigl)
$$
for all $i$, $j$, $k$, $l=1,\cdots, n$ (see equation (16) of \cite{M3}).
We now define
$$
L(R)=T({\cal M}^n(k))/I(R)
$$
which is a bialgebra. $M$ has a right $L(R)$-comodule structure
via the natural projection $T({\cal M}^n(k))\to L(R)$, i.e.
$$
\rho(m_l)=\sum_{v=1}^{n}m_v\ot \overline{c_{vl}}
$$
for all $l=1,\cdots, n$.

First we shall prove the uniqueness of $\sigma$.
Let $\sigma :L(R)\ot L(R)\to k$ such that $(L(R), \sigma)$ is a
Long bialgebra and  $R=R_{\sigma}$. Let $u$, $v=1,\cdots,n$. Then
$$
R_{\sigma}(m_v\ot m_u)=
=\sum_{i,j} \sigma(\overline{c_{iv}}\ot \overline{c_{ju}})m_i\ot m_j
$$
Hence $R_{\sigma}(m_v\ot m_u)=R(m_v\ot m_u)$ gives us
\begin{equation}\label{uns}
\sigma(\overline{c_{iv}}\ot \overline{c_{ju}})=x_{uv}^{ji}
\end{equation}
for all $i$, $j$, $u$, $v=1,\cdots, n$. As $L(R)$ is generated as an
algebra by $(\overline{c_{ij}})$, the relations (\ref{uns}) with
$(L2)-(L5)$ ensure the uniqueness of $\sigma$.

Now we shall prove the existence of $\sigma$. First we define
$$\sigma_0:{\cal M}^n(k)\ot {\cal M}^n(k)\to k,\qquad
\sigma_0(c_{iv}\ot {c_{ju}})=x_{uv}^{ji}
$$
for all $i$, $j$, $u$, $v=1,\cdots, n$.
Then, we extend $\sigma_0$ to a map
$$
\sigma_1:T({\cal M}^n(k))\ot T({\cal M}^n(k))\to k
$$
such that $(L2)-(L5)$ hold. We shall prove that $\sigma_1$ factorizes
to a map
$$\sigma:L(R)\ot L(R)\to k$$
i.e.
$$
\sigma_1\Bigl(T({\cal M}^n(k))\ot o(i,j,k,l) \Bigl)=
\sigma_1\Bigl(o(i,j,k,l)\ot T({\cal M}^n(k))\Bigl)=0.
$$
For $i$, $j$, $k$, $l$, $p$, $q=1,\cdots, n$ we have:
\begin{eqnarray*}
\sigma_1(c_{pq}\ot o(i,j,k,l))&=&
\sum_{v}x_{kv}^{ji}\sigma_1(c_{pq}\ot c_{vl})-
\sum_{\alpha}x_{kl}^{j\alpha}\sigma_1(c_{pq}\ot c_{i\alpha})\\
&=&\sum_{v}x_{kv}^{ji}x_{lq}^{vp}-
\sum_{\alpha}x_{kl}^{j\alpha}x_{\alpha q}^{ip}\\
\text{(from (\ref{sase}) )}
&=&0
\end{eqnarray*}
and
\begin{eqnarray*}
\sigma_1(o(i,j,k,l)\ot c_{pq})&=&
\sum_v x_{kv}^{ji}\sigma_1(c_{vl}\ot c_{pq})-
\sum_{\alpha} x_{kl}^{j\alpha}\sigma_1(c_{i\alpha}\ot c_{pq})\\
&=&\sum_v x_{kv}^{ji}x_{ql}^{pv}-
\sum_{\alpha} x_{kl}^{j\alpha}x_{q\alpha}^{pi} \\
\text{(\mbox{from} (\ref{cinci}) )}
&=&0
\end{eqnarray*}
Hence, we have constructed $\sigma:L(R)\ot L(R)\to k$ such that
$(L2)-(L5)$ hold and $R=R_{\sigma}$. It remains to prove that $\sigma$
satisfies $(L1)$. Using proposition \ref{gen}, it is enough to show that
$(L1)$ holds on the generators. Let $x=\overline{c_{ij}}$,
$y=\overline{c_{pq}}$. We have:
$$
\sum \sigma(x_{(1)}\ot y)x_{(2)}=
\sum_v\sigma(\overline{c_{iv}}\ot \overline{c_{pq}})\:\overline{c_{vj}}=
\sum_vx_{qv}^{pi}\overline{c_{vj}}
$$
and
$$
\sum \sigma(x_{(2)}\ot y)x_{(1)}=
\sum_{\alpha}\sigma(\overline{c_{\alpha j}}\ot \overline{c_{pq}})\:
\overline{c_{i\alpha}}=\sum_{\alpha} x_{qj}^{p\alpha}\overline{c_{i\alpha}}.
$$
Hence,
$$
\sum \sigma(x_{(1)}\ot y)\:x_{(2)}-
\sum \sigma(x_{(2)}\ot y)\:x_{(1)}=\overline{o(i,p,q,j)}=0,
$$
i.e. $(L1)$ holds.
Suppose now that $R$ is bijective and let $S=R^{-1}$. Let $(y_{uv}^{ji})$
be a family of scalars of $k$ such that
$$
S(m_v\ot m_u)=\sum_{i,j}y_{uv}^{ji}m_i\ot m_j,
$$
for all $u$, $v=1,\cdots ,n$.  We define
$$
\sigma^{\prime}:L(R)\ot L(R)\to k, \quad
\sigma^{\prime}(\overline{c_{iv}}\ot \overline{c_{ju}}):=y_{uv}^{ji}
$$
for all $i$, $j$, $u$, $v=1,\cdots, n$. Now it is routine to see that
$\sigma^{\prime}$ is well defined and it is the inverse of $\sigma$.
\end{proof}

\begin{remark}
Let $M$ be a finite dimensional vector space and $R\in \End_k(M\ot M)$
be a solution for the Long equation.
Using theorem \ref{boci} and proposition \ref{bic}, we obtain that
$M$ has a structure of a $L(R)$-dimodule $(M,\cdot,\rho)$ and
$R=R_{(M,\cdot,\rho)}$.
\end{remark}

\begin{examples}
1. Let $a$, $b$, $c\in k$ such that $(b,c)\neq 0$ and $R\in {\cal M}_4(k)$
given by the formula (\ref{100}). Then, $R$ is a solution of the Long
equation. The bialgebra $L(R)$ can be described as follows:

$\bullet$ as an algebra $L(R)=k<x, y>$, the free algebra generated
by $x$ and $y$.

$\bullet$ the comultiplication $\Delta$ and the counity $\varepsilon$
are given by
$$
\Delta(x)=x\ot x, \quad \Delta(y)=x\ot y+y\ot x,\quad
\varepsilon(x)=1, \quad \varepsilon(y)=0.
$$
Indeed, if we write
$$
R(m_v\ot m_u)=\sum_{i,j=1}^{2}x_{uv}^{ji}m_i\ot m_j
$$
we get that, among the elements $(x_{uv}^{ji})$, the only nonzero elements
are:
$$
x_{11}^{11}=ab,\quad x_{21}^{11}=ac,\quad x_{21}^{21}=ab,
\quad x_{12}^{11}=b,
$$
$$
x_{12}^{12}=ab,\quad x_{21}^{11}=c,\quad x_{22}^{21}=b,\quad
x_{22}^{12}=ac,\quad x_{22}^{22}=ab
$$
The sixteen relations $o(i,j,k,l)=0$, written in lexicographic order
according to $(i,j,k,l)$, starting with $(1,1,1,1)$, are
$$
abc_{11}+bc_{21}=abc_{11},\quad abc_{12}+bc_{22}=bc_{11}+abc_{12}
$$
$$
acc_{11}+cc_{21}=acc_{11}, \quad acc_{12}+cc_{22}=cc_{11}+acc_{12}
$$
$$
0=0,\quad 0=0,\quad abc_{11}+bc_{21}=abc_{11},\quad
abc_{12}+bc_{22}=bc_{11}+abc_{12}
$$
$$
abc_{21}=abc_{21},\quad abc_{22}=abc_{22},\quad acc_{21}=acc_{21},\quad
acc_{22}=cc_{21}+acc_{22}
$$
$$
0=0,\quad 0=0,\quad abc_{21}=abc_{21}, \quad abc_{22}=abc_{22}+bc_{21}
$$
Only the four relations below remain:
$$
bc_{21}=0,\quad bc_{22}=bc_{11},\quad
cc_{21}=0,\quad cc_{22}=cc_{11}
$$
As $(b,c)\neq (0,0)$, there are only two linear independent relations:
$$
c_{21}=0,\qquad c_{22}=c_{11}
$$
Now, if we denote $c_{11}=x$, $c_{12}=y$, the above description of
$L(R)$ follows.
Moreover, $L(R)$ can be viewed as a Long bialgebra, where
$\sigma :k<x,y>\ot k<x,y>\to k$ is given by:
$$
\begin{array}{lll}
\sigma(x\ot 1)=1, &\sigma(x\ot x)=ab, & \sigma(x\ot y)=c, \\
\sigma(y\ot 1)=0, &\sigma(y\ot x)=b, & \sigma(y\ot y)=0.
\end{array}
$$
2. Let $a$, $b\in k$ and $R\in {\cal M}_4(k)$ given by
$$
R=
\left(
\begin{array}{cccc}
a&0&0&0\\
0&b&0&0\\
0&0&0&0\\
0&0&0&0
\end{array}
\right)
$$
Then $R$ is a solution of the Long equation, as $R=f\ot g$,
where $f$ and $g$ are given by
$$
f=
\left(
\begin{array}{cc}
1&0\\
0&0
\end{array}
\right)
\qquad
g=
\left(
\begin{array}{cc}
a&0\\
0&b
\end{array}
\right)
$$
i.e. $fg=gf$. We shall describe the bialgebra $L(R)$.
Suppose $(a,b)\neq 0$ (otherwise
$R=0$ and $L(R)=T({\cal M}^4(k))$ ).
Among the sixteen relations $o(i,j,k,l)=0$, the only linear independent
ones are:
$$
ac_{12}=ac_{21}=bc_{12}=bc_{21}=0
$$
As $(a,b)\neq 0$, we obtain $c_{12}=c_{21}=0$. Now, if we denote
$c_{11}=x$, $c_{22}=y$, we get that the bialgebra $L(R)$ has the
following description:

$\bullet$ as an algebra, $L(R)=k<x, y>$, the free algebra generated
by $x$ and $y$.

$\bullet$ the comultiplication $\Delta$ and the counity $\varepsilon$
are given in such a way that $x$ and $y$ are groupal elements
$$
\Delta(x)=x\ot x, \quad \Delta(y)=y\ot y\quad
\varepsilon(x)=\varepsilon(y)=1.
$$
Moreover, $k<x,y>$ is a Long bialgebra as follows:
$$
\begin{array}{lll}
\sigma(x\ot 1)=1, &\sigma(x\ot x)=a, & \sigma(x\ot y)=b, \\
\sigma(y\ot 1)=1, &\sigma(y\ot x)=0 & \sigma(y\ot y)=0.
\end{array}
$$
3. Let $M$ be an $n$-dimensional vector space with $\{m_1,\cdots, m_n\}$
a basis and $\phi:\{1,\cdots,n \}\to \{1,\cdots,n \}$ a function
with $\phi^{2}=\phi$. We have proven that the map
\begin{equation}\label{dor}
R^{\phi}:M\ot M\to M\ot M,\quad
R^{\phi}(m_i\ot m_j)=\delta_{ij}\delta_{i\in Im(\phi)}
\sum_{a,b\in \phi^{-1}(i)}m_a\ot m_b
\end{equation}
for all $i$, $j=1,\cdots, n$, is a solution of the Long equation.
Below, we shall describle the bialgebra $L(R^{\phi})$. From equation
$(\ref{dor})$, the scalars $(x_{uv}^{ji})$ which define $R^{\phi}$ are given
by
$$
x_{uv}^{ji}=\delta_{uv}\delta_{\phi(i)v}\delta_{\phi(j)v}
$$
for all $i$, $j$, $u$, $v=1,\cdots ,n$.
The relations $o(i,j,k,l)=0$ are
$$
\sum_{v}\delta_{kv}\delta_{\phi(i)v}\delta_{\phi(j)v}c_{vl}
=\delta_{kl}\delta_{\phi(j)l}\sum_{\alpha}\delta_{\phi(\alpha)l}c_{i\alpha}
$$
or equivalent to
$$
\delta_{\phi(i)k}\delta_{\phi(j)k}c_{kl}=
\delta_{kl}\delta_{\phi(j)l}\sum_{\alpha}\delta_{\phi(\alpha)l}c_{i\alpha}
$$
for all $i$, $j$, $u$, $v=1,\cdots ,n$.
All the relations $o(i,j,k,l)=0$ are identities $0=0$, with the exception
of the following three types of relations:
\begin{equation}
\left\lbrace
\begin{array}{llll}
o(i,j,\phi(j), l)=0,\quad \mbox{for all}\;\;\; l\neq \phi(j),\; \phi(i)=\phi(j)\\
o(i,j,\phi(j), \phi(j))=0,\quad \mbox{for all}\;\;\; \phi(i)\neq\phi(j)\\
o(i,j,\phi(j), \phi(j))=0,\quad \mbox{for all}\;\;\; \phi(i)=\phi(j)
\end{array}
\right.
\end{equation}
which give us the following description of the bialgebra $L(R^{\phi})$: as
an algebra, $L(R^{\phi})$ is the free algebra generated by $(c_{ij})$ with
the relations:

\begin{equation}\label{rel}
\left\lbrace
\begin{array}{llll}
c_{\phi(j)l}=0, \quad \mbox{for all}\;\;\;  l\neq \phi(j)\\
\sum_{\alpha\in \phi^{-1}(\phi(j))}c_{i\alpha}=0, \quad \mbox{for all}\;\;\;
\phi(i)\neq\phi(j)\\
\sum_{\alpha\in \phi^{-1}(\phi(i))}c_{i\alpha}=
c_{\phi(i)\phi(i)},\quad \mbox{for all}\;\;\;  i=1,\cdots, n
\end{array}
\right.
\end{equation}

The comultiplication and the counity of $L(R^{\phi})$ are given in such
a way that the matrix $(c_{ij})$ is comultiplicative.
The bialgebra $L(R^{\phi})$ is a Long bialgebra with
$$
\sigma:L(R^{\phi})\ot L(R^{\phi}) \to k, \quad
\sigma(c_{iv}\ot c_{ju})=\delta_{uv}\delta_{\phi(i)v}\delta_{\phi(j)v}
$$
for all $i$, $j$, $u$, $v=1,\cdots ,n$.

4. In the next three examples we point out particular cases of the
previous example. Let $\phi$ the identity map, i.e. $\phi(i)=i$ for all
$i=1,\cdots, n$. Then $L(R^{\phi})=k<x_1,\cdots, x_n>$, the free algebra
generated by $x_1,\cdots, x_n$, and the coalgebra structure is given in
such a way that $x_i$ is a group-like element, i.e.
$\Delta(x_i)=x_i\ot x_i$ and $\varepsilon(x_i)=1$ for all $i=1,\cdots, n$.
Indeed, from the  relations (\ref{rel}), only
$$
c_{jl}=0, \forall l\neq j
$$
remains. If we denote $c_{ii}=x_i$ the conclusion follows.
We note that the bialgebra $k<x_1,\cdots, x_n>$ is a
Long bialgebra with $\sigma: k<x_1,\cdots, x_n>\ot k<x_1,\cdots, x_n>\to k$
given by
$$
\sigma(x_i\ot x_j)=\delta_{ij}
$$
for all $i$, $j=1,\cdots, n$.

5. Let $n=4$ and $\phi$ given by
$$
\phi(1)=1,\quad \phi(2)=\phi(3)=\phi(4)=2.
$$
Then the relations (\ref{rel}) are:
$$
c_{1l}=c_{2j}=c_{i1}=0, \quad
c_{32} + c_{33} + c_{34}= c_{22},\quad
c_{42} + c_{43} + c_{44}= c_{22}
$$
for all $l\neq 1$, $j\neq 2$, $i\neq 1$. Now, if we denote
$c_{11}=x_{1}$, $c_{22}=x_{2}$, $c_{32}=x_{3}$,
$c_{33}=x_{4}$, $c_{42}=x_{5}$, $c_{44}=x_{6}$ we obtain the description
of the corresponding bialgebra $L(R^{\phi})$:

$\bullet$ as an algebra $L(R^{\phi})=k<x_{i}\mid i=1,\cdots, 6 >$,
the free algebra generated by six generators.

$\bullet$ the comultiplication $\Delta$ and the counity $\varepsilon$
are given by
$$
\Delta (x_{1})=x_{1}\ot x_{1},\quad
\Delta (x_{2})=x_{2}\ot x_{2},
$$
$$
\Delta (x_{3})=x_{3}\ot x_{2}+x_{4}\ot x_{3}+ (x_2-x_3-x_4)\ot x_{5},\quad
\Delta (x_{4})=x_{4}\ot x_{4}+(x_2-x_3-x_4)\ot (x_2-x_5-x_6),
$$
$$
\Delta (x_{5})=x_{5}\ot x_{2}+ (x_2-x_5-x_6)\ot x_3+x_6\ot x_5,\quad
\Delta (x_{6})=(x_2-x_5-x_6)\ot (x_2-x_3-x_4)+x_{6}\ot x_{6}
$$
$$
\varepsilon(x_{1})=\varepsilon(x_{2})=\varepsilon(x_{4})=
\varepsilon(x_{6})=1,\quad
\varepsilon(x_{3})=\varepsilon(x_{5})=0.
$$
6. Let $n=4$ and $\phi$ given by
$$
\phi(1)=\phi(2)=2,\quad \phi(3)=\phi(4)=4.
$$
Then the relations (\ref{rel}) are:
$$
c_{2l}=c_{4j}=0,\quad
c_{11}+c_{12}=c_{22},\quad
c_{33}+c_{34}=c_{44},\quad
c_{13}+c_{14}=c_{31}+c_{32}=0
$$
for all $l\neq 2$, $j\neq 4$. If we denote
$c_{11}=x_{1}$, $c_{12}=x_{2}$, $c_{13}=x_{3}$,
$c_{31}=x_{4}$, $c_{33}=x_{5}$, $c_{34}=x_{6}$, we obtain the description
of the corresponding bialgebra $L(R^{\phi})$:

$\bullet$ as an algebra $L(R^{\phi})=k<x_{i}\mid i=1,\cdots, 6 >$,
the free algebra generated by six generators.

$\bullet$ the comultiplication $\Delta$ and the counity $\varepsilon$
are given by
$$
\Delta (x_{1})=x_{1}\ot x_{1}+x_3\ot x_4,\quad
\Delta (x_{2})=x_{1}\ot x_{2}+x_2\ot x_1+x_2\ot x_2-x_3\ot x_4,
$$
$$
\Delta (x_{3})=x_{1}\ot x_{3}+x_{3}\ot x_{5},\quad
\Delta (x_{4})=x_{4}\ot x_{1}+x_5\ot x_4,
$$
$$
\Delta (x_{5})=x_{4}\ot x_{3}+ x_5\ot x_5,\quad
\Delta (x_{6})=-x_4\ot x_3+x_5\ot x_6+x_6\ot x_5+x_{6}\ot x_{6}
$$
$$
\varepsilon(x_{1})=\varepsilon(x_{5})=1,\quad
\varepsilon(x_{2})=\varepsilon(x_{3})=
\varepsilon(x_{4})=\varepsilon(x_6)=0.
$$
\end{examples}

\end{document}